\documentclass [11pt,titlepage]{article}

\usepackage{amsmath}
\usepackage{amssymb}
\begin{document}

\begin{center}
\begin{Large} 
{\bf The diophantine equation $\frac{x^3}{3}+y^3+z^3-2xyz=0$}  

\end{Large}

\addvspace{0.5cm}
{\bf Joseph Amal Nathan}

Reactor Physics Design Division, Bhabha Atomic Research Centre, Mumbai-400085, India\\
{\bf email:}josephan@magnum.barc.ernet.in
\end{center}
\noindent
{\bf Abstract:}~We will be presenting two theorems in this paper. The first theorem, which is a new result, is about the non-existence of integer solutions of the cubic diophantine equation. In the proof of this theorem we have used some known results from theory of binary cubic forms and the method of infinite descent, which are well understood in the purview of  Elementary Number Theory(ENT). In the second theorem, we show, that the famous Fermat's Last Theorem(FLT) for exponent $3$ and the first theorem are equivalent. So Theorem1 and 2 constitute an alternate proof for the non-existence of integer solutions of this famous cubic Fermat's equation. It is well known, that from L.Euler($1770$) to F.J.Duarte($1944$) many had given proof of FLT for exponent $3$. But all proofs uses concepts, which are beyond the scope of ENT. Hence unlike other proofs the proof given here is as an ENT proof of FLT for exponent $3$.

\noindent
We will denote the greatest common divisor of integers $a_1,~a_2,~a_3 ,...$ by symbol $(a_1,~a_2,~a_3 ,...)$. For the results from theory of Binary Cubic Forms, we follow L.J.Mordell[1]. Consider the binary cubic
$$f(x,y)=a x^3+b x^2 y+c x y^2+d y^3= \{ a,b,c,d \} ,$$
\noindent
with integer coefficients and discriminant 
$$D=-27 a^2 d^2 +18 a b c d +b^2 c^2 - 4 a c^3 -4 b^3 d,$$
\noindent
where $D \neq 0$. The quadratic covariant $H(x,y)$ is
\begin{eqnarray}
\nonumber H(x,y) &\!\!\!=\!\!\!& (b^2-3 a c)x^2+(b c-9 a d) x y+(c^2-3 b d)y^2,\\
&\!\!\!=\!\!\!& Ax^2+Bxy+Cy^2=\{ A,B,C \}
\end{eqnarray}
\noindent
with discriminant $B^2-4AC=-3D$ and the cubic covariant
\begin{eqnarray*}
G(x,y) &\!\!\!=\!\!\!& -(27a^2d-9abc+2b^3)x^3+3(6ac^2+b^2c-9abd)x^2y\\
&&~~~~~+3(bc^2-6b^2d+9acd)xy^2+(27ad^2-9bcd+2c^3)y^3.
\end{eqnarray*}
\noindent
{\bf Lemma1:} $f(x,y),H(x,y)$ and $G(x,y)$ {\it are algebraically related by the identity}, $$G^2(x,y)+27Df^2(x,y)=4H^3(x,y)~\hbox{[1]}.$$
%\newpage
\noindent
{\bf Lemma2:} {\it All integer solutions of}
\begin{equation}
X^2+27 k  Y^2=4Z^3,~~~~~~~~~~(X,Z)=1
\end{equation}
\noindent
{\it are given by taking,} $X=G(x,y),~k =D,~Y=f(x,y),~Z=H(x,y).$

\noindent
{\it Proof}.~For proof of this lemma, we follow word by word the proof given in L.J.Mordell[1] for equation of the form $X^2+k Y^2=Z^3$.

\noindent
Let $[X,Y,Z]=[g,f,h]$ be a solution such that
\begin{equation}
g^2+27k f^2=4h^3,~~~~~~~~~~(g,h)=1.
\end{equation}
\noindent
We will construct a binary cubic $f(x,y)$ with integer coefficients and of discriminant $D=k $ such that $g,~f,~h$ are values assumed by $G(x,y),~f(x,y),~H(x,y)$ for integer $x,~y$. Then all solutions of (2) are given by taking $f(x,y)$ a set of binary cubics of discriminant $D$ and letting $x,~y$ run through all integer values for which $(X,Z)=1$.

\noindent
Since from (3), $-3k $ is a quadratic residue of $h$, there exist binary quadratics of discriminant $-3D$ with first coefficient $h$,
$$ \{ h,B,C \} =hx^2+Bxy+Cy^2,~~~~~\hbox{where}~~B^2-4hc=-3k .$$
\noindent
We take for $B$ any solution of the congruence $3fB \equiv -g(\bmod~4h^3)$. We will now construct a binary cubic $ \{ f,b,c,d \} $ with discriminant $D=k $ and $H(x,y)$ given in (1) by $\{ h,B,C \} $. Since,
$$h=b^2-3fc,~~~~~~~~~~c=\frac{b^2-h}{3f},$$
\noindent
we can take $b \equiv g/2h(\bmod~f)$ and in particular, $2bh=g+3fB$, and so $b \equiv 0(\bmod~2h^2)$. Then
\begin{eqnarray*}
4h^2c&\!\!\!=\!\!\!&\frac{(g+3fB)^2-4h^3}{3f}\\
&\!\!\!=\!\!\!&-9k f+2gB+3fb^2.
\end{eqnarray*}
\noindent
We now show $c$ is an integer,
\begin{eqnarray*}
3f~4h^2c &\!\!\!\equiv\!\!\!& 3f\left[ -9k f+2g \left( -\frac{g}{3f} \right)+3f \left( \frac{g}{3f}\right)^2 \right] (\bmod~h^2)\\
&\!\!\!\equiv\!\!\!& -27k f^2-g^2 (\bmod~h^2) \equiv 0 (\bmod~h^2).
\end{eqnarray*}
\noindent
We now find $d$. Since $bc-9fd=B$,
$$9fd=\left(\frac{g+3fB}{2h}\right)\left(\frac{-9k f+2gB+3fb^2}{4h^2}\right)-B.$$ 
\noindent
We will simplify the above equation and show $d$ is an integer for completion of proof.
$$24h^3d=-3kg-27k fB+3gB^2+3fB^3,$$
\begin{eqnarray*}
~~~~~~~~~~~~~~~~~~~~~~~~~9f^2~24h^3d&\!\!\!\equiv\!\!\!&9f^2 \left[ -3k g+9k g+\frac{3g^3}{9f^2}-\frac{g^3}{9f^2} \right](\bmod~h^3)\\
&\!\!\!\equiv\!\!\!& 2g \left[ 27k f^2+g^2 \right] (\bmod~h^3) \equiv 0(\bmod~h^3).~~~~~~~~~~~~~~~~~~~~~~~~~~~\Box
\end{eqnarray*}
{\bf Lemma3:} {\it For} $D=1$ {\it the classes of binary cubics are given by} $a=0,~|b|= 1,~|c|=1,~d=0$.

\noindent
{\it Proof}.~The discriminant $D$ is also the invariant
$$D=a^4 (\alpha - \beta)^2 (\beta - \gamma)^2 (\gamma - \alpha)^2,$$
\noindent
where $\alpha ,~\beta ,~\gamma$ are the roots of the equation
$$a \zeta^3 +b \zeta +c \zeta +d=0.$$
\noindent
When $D>0$ the roots $\alpha ,~\beta ,~\gamma$ are all real and distinct. So from (1) for expressions of $A$ we get
%\begin{eqnarray}
$$\frac{A}{a^2}=\left( \frac{b}{a} \right)^2-3 \left( \frac{c}{a} \right)=(\alpha ^2+\beta ^2+\gamma ^2)-(\alpha \beta +\beta \gamma +\gamma \alpha ) > 0$$
%\left( \frac{c}{a} \right)^2 -3 \left( \frac{b}{a} \right) \left( \frac{d}{a} \right) %&=&(\alpha ^2 \beta ^2+\beta ^2 \gamma ^2+\gamma ^2 \alpha ^2)-(\alpha \beta \gamma )(\alpha %+\beta +\gamma ) > 0.
%\end{eqnarray*}
\noindent
Since $D,~A>0$, if $H(x,y)=\frac{1}{4A} \left[ (2Ax+By)^2+3Dy^2 \right]=0 \Rightarrow x=y=0$. Hence for $D>0$ we have, $B^2-4AC<0$ and $A>0$. Then $H(x,y)$ is positive definite. We know a binary quadratic form $\{ I,J,K \} $ is reduced if $K \geq I \geq |J|$. Since every class of positive definite binary forms contains at least one reduced form, we can take $C \geq A \geq |B|$. So $AC \geq B^2$, from $4AC-B^2=3D$, we get $AC \leq D$ and $A \leq \sqrt D$. 

\noindent
Now for $D=1$ the above inequalities give, $A \leq 1,~C \leq 1$ and $A,~C\geq |B| \Rightarrow |B| \leq 1.$ Since $A,~C>0$, $A=C=1$. From $4AC-B^2=3D$ we get $B^2=1 \Rightarrow |B|=|bc-9ad|=1$, so $|b|,~|c| \geq 1$, since none of $b,~c$ can be zero. Substituting for $a,~d$ from $b^2-3ac=1,~c^2-3bd=1$, 
in $|B|=1$ we get, 
$$b^2+c^2-1=\pm bc.$$
\noindent
Let $|b| \geq |c|>1$. From above equation $b^2+c^2-1 \leq |bc| \Rightarrow b^2+c^2-1 \leq b^2 \Rightarrow c^2 \leq 1$ contradicting $|c|>1$. So $|b|=|c|=1$ and $a=d=0$.~~~~~~~~~~~~~~~~~~~~~~~~~~~~~~~~~~~~~~~~~~~~~~~~~~~~~~~~~~~~$\Box$

\noindent
Later in this paper, we will require (2) with $k=1$. From Lemma3 for $D=1$ we have $|b|=|c|=1$ and $a=d=0$. Since $k=D=1$ from the form of (2), we see, it is sufficient to take $a=0,~b=\pm 1,~c=1,~d=0$. So we have,
\begin{equation}
f(x,y)=xy(x\pm y),~~~~~H(x,y)=x^2\pm xy+y^2,~~~~~G(x,y)=2x^3\pm 3x^2y-3xy^2\mp 2y^3.
\end{equation}
\noindent
{\bf Lemma4:} {\it Given nonzero integers $x,~y$ such that $(x,y)=1$ for any positive odd integer $n$},
$$\left( x+y, \frac {x^n+y^n}{x+y} \right) = (x+y,n).$$ 
\noindent
{\it Proof}.~Let $(x+y,n)=g$. Since, $\frac {x^n+y^n}{g(x+y)} = \frac{(x+y)}{g} \left[\sum_{i=1}^{n-1}(-1)^{i-1}~i~x^{n-i-1} y^{i-1} \right]+\frac{n}{g}y^{n-1},$ any divisor of $\frac{(x+y)}{g}$ will divide all terms except last. Hence $\left( \frac{x+y}{g}, \frac {x^n+y^n}{g(x+y)} \right) = 1.~~~~~~~~~~~~~~~~~~~~~~~~~~~~~~~~~~~~~~~~~~~~~~~~~~~~\Box$
\newpage
\noindent
{\bf Theorem1:} {\it The cubic equation}
\begin{equation}
\frac{u^3}{3}+v^3+w^3-2uvw=0, 
\end{equation}
{\it has no solution in nonzero integers} $u,~v,~w$.

\noindent
{\it Proof}:~There is no loss of generality if $u,~v,~w$ are pairwise prime integers. From (5) we see $3 \mid u$ and $u \mid (v^3+w^3) \Rightarrow 9\mid (v^3+w^3) \Rightarrow 9 \mid u$. 

\noindent
We can always choose nonzero integers $M,~N,~p,~q$ such that,
$$v+w=3Mp,~~v^2-vw+w^2=9M^2p^2-3vw=3Nq,~~u=9pq,$$
\noindent
where $3 \nmid M$ and $3 \mid p$ if and only if $3^3 \mid u$. One can see $N,~q$ are odd integers. By Lemma4 we have $3 \parallel (v^2-vw+w^2)$ hence,
\begin{equation}
3 \nmid N,~q~~\hbox{and}~~(M,N)=(M,q)=(p,N)=(p,q)=1.
\end{equation}
\noindent
Substituting for $u,~v+w,~vw$ in terms of $M,~N,~p,~q$ in (5) and after rearranging we get
\begin{equation}
3p^2(9q^2-2M^2)=-N(M+2q).
\end{equation}
\noindent
Now using (7) we will find $M,~N$ in terms of $p,~q$.
Since $q$ is odd $9q^2-2M^2$ is odd. Let $(9q^2-2M^2,M+2q)=\delta $. Since $9q^2-2M^2=q^2+2(4q^2-M^2)$ any prime divisor $\epsilon $ of $\delta $ divides $q \Rightarrow \epsilon \mid M$. But $(M,q)=1$ hence $\delta =1$. Since $(N,3p^2)=1$ from (7) we have,
$$M=3p^2-2q~~~~~\hbox{and}~~~~~N=2M^2-9q^2=18p^4-24p^2q-q^2.$$
\noindent
We can see $u,~v+w,~vw$ are integers. If $v-w$ is also an integer it shows the existence of integer solutions for (5) in terms of integers $p,q$. Substituting for $v+w,~vw$ in the following expression we get,
$$(v-w)^2=(v+w)^2-4vw=-4q^3-27[p(p^2-2q)]^2,$$
\noindent
which can be written as,
$$(v-w)^2+27[p(p^2+2(-q))]^2=4(-q)^3.$$
\noindent
From Lemma2 all integer solutions of the above equation is given by taking,
$$X=G(x,y)=v-w,~~~k=D=1,~~~Y=f(x,y)=p(p^2-2q),~~~Z=H(x,y)=-q.$$
\noindent
Since $k=D=1$, substituting for $f(x,y),~H(x,y)$ from (4) in above equation for $Y$,
\begin{equation}
p[p^2+2(x^2\pm xy+y^2)]=xy(x\pm y)
\end{equation}
\noindent
Let $(x,y)=\delta $. $q=-(x^2 \pm xy+y^2)$ so any prime divisor $\epsilon $ of $\delta $ divides $q$. From (8) we see $\epsilon $ divides either $p$ or $p^2-2q \Rightarrow \epsilon \mid p$. But $(p,q)=1$ hence $\delta =1$. 

\noindent
To find the integer solutions of (8), we follow a similar procedure as we did for (5). It is possible to choose nonzero integers $Q,~R,~S,~T$ such that 
$$QR=xy,~~~~~ST=x\pm y,~~~~~p=RT.$$
\noindent
Since $(x,y)=1$, we have $(Q,S)=(Q,T)=(R,S)=(R,T)=1$. Substituting for $x,~y$ in terms of $Q,~R,~S,~T$ in (8) and after simplifying we get,
\begin{equation}
T^2(R^2+2S^2)=Q(S \pm 2R).
\end{equation}
\noindent
Now we will show $3 \mid f(x,y)$. From (6), we have $3 \nmid q$. Assume $3 \nmid f(x,y)$ then ${x \equiv \pm y~(\bmod~3)}$ $ \Rightarrow q \equiv 0~(\bmod~3)$ a contradiction. Now let $(Q,R)=\delta $ any prime divisor $\epsilon $ of $\delta $  from (9) will divide $(R^2+2S^2)$ since $(Q,T)=1$. Hence $\epsilon $ has to divide $2$ or $S$. Since $(Q,S)=1$ we have $\epsilon \leq 2 \Rightarrow (Q,R) \leq 2$. Similarly let $(S,T)=\delta $ any prime divisor $\epsilon $ of $\delta $  again from (9) will divide $(S \pm 2R)$ since $(Q,T)=1$. So $\epsilon $ has to divide $2$ or $R$. Since $(R,T)=1$ we get $(S,T) \leq 2$. So we see $3 \mid Q$ or $R$ or $S$ or $T$.
 
\noindent
To evaluate $Q,~S$ in terms of $R,~T$ using (9) we should know about $(R^2+2S^2,S \pm 2R)$. Let $(R^2+2S^2,S \pm 2R)=\Delta$ and $\eta $ be a prime divisor of $\Delta$. Since $R^2+2S^2=2(S-2R)(S+2R)+9R^2$, $\eta \mid 3$ or $R$. If $\eta \mid R$ then $\eta \mid S \Rightarrow \Delta=1 $ since $(R,S)=1$. If $\eta \mid 3$ then $\Delta $ can only be a power of $3$.
 
\noindent
Case-I: For following cases $\Delta =1$.

\noindent
(i) $3 \mid Q,~(S \pm 2R)$. From (9) we choose the case $T^2[2(S-2R)(S+2R)+9R^2]=Q(S \mp 2R)$. Since $(R,S)=1$ we have $(2R,S) \leq 2$ and so $3 \nmid (S \mp 2R)$. Any prime divisor $\eta $ of $\Delta $ divides $R$ and from above argument we get $\Delta =1$. The other part $T^2(R^2+2S^2)=Q(S \pm 2R)$ when $3 \mid Q,~(S \pm 2R)$ will be taken in Case-II. 

\noindent
(ii) $3 \mid R$ or $S$. In both the cases $3 \nmid (S \pm 2R)$. Since $R^2+2S^2=2(S^2-4R^2)+9R^2$ any prime divisor of $\Delta $ divides $R$ so $\Delta =1$.

\noindent
For (i) and (ii) we have $(R^2+2S^2,S \pm 2R)=1$ and $(Q,T)=1$. So from (9)
$$S=T^2 \mp 2R,~~~~~Q=R^2+2S^2=9R^2+2T^4 \mp 8 T^2 R.$$
\noindent
As we did before when we evaluate $x \mp y$ in terms of $R,T$ we get,
\begin{equation}
(x \mp y)^2=(x \pm y)^2 \mp 4xy=\mp 36 R^3+T^2(T^2 \mp 6R)^2.
\end{equation}
\noindent
We have $(R,T)=1$ and $3 \nmid T$. Let $(36 R^3,T^2(T^2 \mp 6R)^2)=\delta $. Let $R$ be even or odd and $T$ be odd then $T(T^2 \mp 6PL)$ is odd. From (10) $(x \mp y)$ is odd and $2 \nmid \delta$, then any prime divisor of $\delta $ has to divide $T \Rightarrow \delta =1$. If $R$ is odd and $T$ is even, $4 \mid T(T^2 \mp 6R)$. From (10) $2 \| (x \mp y)$ and $\delta =4$. Since $\delta =1$ or $4 \Rightarrow \left( (x \mp y)^2,T^2(T^2 \mp 6R)^2 \right)=1$ or $4$ respectively we get $\left( (x \mp y)+T(T^2 \mp 6PL),(x \mp y)-T(T^2 \mp 6PL) \right)=2$. We will see there will be no loss of generality if we take $(x \mp y)-T(T^2 \mp 6R)$ is divisible by $3$. Now from (10) we get the following equations,
\begin{equation}
\begin{array}{ccccccc}
(x-y)+T(T^2-6PL)&\!\!\!=\!\!\!&\pm 2P^3&~~~\hbox{and}~~~&(x-y)-T(T^2-6PL)&\!\!\! =\!\!\! &\mp 18L^3,
\end{array}
\end{equation}
\begin{equation}
\begin{array}{ccccccc}
(x+y)+T(T^2+6PL)&\!\!\!=\!\!\!&+2P^3&~~~\hbox{and}~~~&(x+y)-T(T^2+6PL)&\!\!\! =\!\!\! &+18L^3,
\end{array}
\end{equation}
\noindent
where $R=PL$ and $(P,L)=1$. Eliminating $(x-y)$ from (11) and $(x+y)$ from (12),
$$T(T^2-6PL)= \pm P^3 \pm 9L^3,~~~~~T(T^2+6PL)=P^3-9L^3.$$
\noindent
Substituting for $\{ 3L= \mp r,P= \mp s,T=t \}$ and $\{ 3L=r,P=-s,T=t \}$ in the above first and second equations respectively we get,
$$\frac{r^3}{3}+s^3+t^3-2rst=0,$$
\noindent
which has the form of (5). Since we have $u=9pq$, $p=RT$, $R=PL$, $|3L|=|r|$, we see
$$|u|>|3p|\geq|3R|\geq|3L|=|r|.$$
\noindent
Hence by the method of infinite descent, there will be no solution in nonzero integers $u,~v,~w$ for (5).

\noindent
Now we take the remaining cases for the completion of the proof.

\noindent
Case-II: For remaining cases $\Delta =3^2$.

\noindent
(i) $3 \mid Q,~(S \pm 2R)$. We have to take $T^2[2(S-2R)(S+2R)+9R^2]=Q(S \pm 2R)$. We have $3 \nmid (S \mp 2R)$ and $3^2 \mid (S \pm 2R)$. If $3^3 \mid (S \pm 2R)$ then $3^4 \mid Q(S \pm 2R)$ but $3^2 \| (R^2+2S^2) \Rightarrow 3 \mid T$ a contradiction. So $3^2 \| (S \pm 2R)$ which gives $\Delta =3^2$. 

\noindent
(ii) $3 \mid T$. Since $(Q,T)=1$ we see from (9) $3^2 \mid (S \pm 2R),(R^2+2S^2)$ and again from (9) we get $3^4 \mid (S \pm 2R) \Rightarrow 3^2\| (R^2+2S^2) \Rightarrow \Delta =3^2$ 

\noindent
For (i) and (ii) we have $(R^2+2S^2,S \pm 2R)=3^2$ and $(Q,T)=1$. Hence from (9) we get
$$S=9T^2 \mp 2R,~~~~~9Q=R^2+2S^2=R^2+18T^4 \mp 8 T^2 R.$$
\noindent
Evaluating $x \mp y$ in terms of $R,T$ we get,
\begin{equation}
(x \mp y)^2=\mp 4R^3+T^2(9T^2 \mp 6R)^2.
\end{equation}
\noindent
We have $(R,T)=1$ and $3 \nmid R$. Let $(4R^3,T^2(9T^2 \mp 6R)^2)=\delta $. Let $R$ be even or odd and $T$ be odd then $T(9T^2 \mp 6PL)$ is odd. From (13) $(x \mp y)$ is odd and $2 \nmid \delta$, then any prime divisor of $\delta $ has to divide $T \Rightarrow \delta =1$. If $R$ is odd and $T$ is even, $4 \mid T(9T^2 \mp 6R)$. From (13) $2 \| (x \mp y)$ and $\delta =4$. Since $\delta =1$ or $4 \Rightarrow \left( (x \mp y)^2,T^2(9T^2 \mp 6R)^2 \right)=1$ or $4$ respectively we get $\left( (x \mp y)+T(9T^2 \mp 6PL),(x \mp y)-T(9T^2 \mp 6PL) \right)=2$. From (13), 
\begin{equation}
\begin{array}{ccccccc}
(x-y)+T(9T^2-6PL)&\!\!\!=\!\!\!& \pm 2P^3&~~~\hbox{and}~~~&(x-y)-T(9T^2-6PL)&\!\!\!=\!\!\!& \mp 2L^3,
\end{array}
\end{equation}
\begin{equation}
\begin{array}{ccccccc}
(x+y)+T(9T^2+6PL)&\!\!\!=\!\!\!&+2P^3&~~~\hbox{and}~~~&(x+y)-T(9T^2+6PL)&\!\!\!=\!\!\!&+2L^3,
\end{array}
\end{equation}
\noindent
where $R=PL$ and $(P,L)=1$. As before eliminating $(x-y)$ from (14) and $(x+y)$ from (15),
$$T(9T^2-6PL)= \pm P^3 \pm L^3,~~~~~T(9T^2+6PL)=P^3-L^3.$$
\noindent
Substituting for $\{ 3T=r,P= \mp s,L= \mp t \}$ and $\{ 3T=r,P=-s,L=t \}$ in above first and second equations respectively we again arrive at,
$$\frac{r^3}{3}+s^3+t^3-2rst=0,$$
\noindent
We have $u=9pq$, $p=RT$, $|3T|=|r|$, so
$$|u|>|3p|\geq|3T|=|r|.$$
\noindent
Again there are no solution in nonzero integers $u,~v,~w$ for (5).~~~~~~~~~~~~~~~~~~~~~~~~~~~~~~~~~~~~~~~~~~~$\Box$

\noindent
R. Perrin[2] had shown the following fact concerning FLT for exponent $3$,

\noindent
{\it The following statements are equivalent and true:}

\noindent
(1) {\it Fermat's last theorem is true for the exponent} $3$.\\
\noindent
(2) {\it For every} $n \geq 1$ {\it the equation}
$$X^3+Y^3+3^{3n-1}Z^3=2~3^n X Y Z$$
\noindent
{\it has no solution in nonzero integers $X,~Y,~Z$, not multiples of} $3$.

\noindent
We will show that the above equation and (5) are same. In (5) we have $3 \mid u$ and $3 \nmid v,~w$. Let $n \geq 1$ such that $3^n \| u$. If we take $u=3^nZ,~v=X,~w=Y$ we see $3 \nmid XYZ$. After substituting for $u,~v,~w$ in (5) we get the above equation. 

\noindent
Now we state the following theorem wherein we show Theorem1 and FLT for exponent $3$ are equivalent by an independent method.

\noindent
{\bf Theorem2:} {\it The equation} $x^3+y^3+z^3=0$ {\it has no solution in nonzero integers} $x,~y,~z$.

\noindent
{\it Proof}:~There is no loss in generality if $x,~y,~z$ are pairwise prime. 
Define $m=x+y+z$, we see $3 \mid m$. From $[m-(x+y)]^3=z^3$ we get
$$m^3-3m^2(x+y)+3m(x+y)^2-3xy(x+y)=0 \Rightarrow 3 \mid xyz.$$
\noindent
Let $3 \mid z$. Now $-z^3=(x+y)(x^2-xy+y^2)$, from Lemma4,
\begin{eqnarray*}
\begin{array}{ccccccccc}
3(x+y)&\!\!\!=\!\!\!&u^3,~~~~~~~~~~& x^2-xy+y^2&\!\!\!=\!\!\!&3U^3,~~~~~~~~~~& z&\!\!\!=\!\!\!&-uU  \\
~(z+x) &\!\!\!=\!\!\!&v^3,~~~~~~~~~~& z^2-zx+x^2&\!\!\!=\!\!\!&V^3,~~~~~~~~~~& y&\!\!\!=\!\!\!&-vV   \\
~(z+y) &\!\!\!=\!\!\!&w^3,~~~~~~~~~~& y^2-zy+z^2&\!\!\!=\!\!\!&W^3,~~~~~~~~~~& x&\!\!\!=\!\!\!&-wW
\end{array}
\end{eqnarray*}
\noindent
where $u,~v,~w,~U,~V,~W$ are relatively prime nonzero integers and $3 \mid u$. From the above equations we get,
\begin{equation}
2z=- \frac{u^3}{3}+v^3+w^3,~~2y= \frac{u^3}{3}-v^3+w^3,~~2x= \frac{u^3}{3}+v^3-w^3,~~2m= \frac{u^3}{3}+v^3+w^3.
\end{equation}
\noindent
Also from the definition of $m$ we get
\begin{equation}
m=u \left( \frac{u^2}{3}-U \right) =v(v^2-V)=w(w^2-W) \Rightarrow m= \Gamma uvw,
\end{equation}
\noindent
where $\Gamma= \left( \frac{u^2}{3}-U,v^2-V,w^2-W \right)$. To evaluate $\Gamma$, after substituting for $z,~y,~x,~m$, from (16), in the following expression, we get,
$$\Gamma=\sqrt[3]{\frac{x^3+y^3+z^3+m^3}{u^3 v^3 w^3}}=1.$$
\noindent
Equating the expressions of $m$ in (16) and (17) we get,
$$\frac{u^3}{3}+v^3+w^3-2uvw=0,$$ 
\noindent
where $u,~v,~w$ are pairwise relatively prime integers. Now from Theorem1 there are no nonzero integers $u,~v,~w$ satifying the above equation. Hence the equation $x^3+y^3+z^3=0$ has only trivial solutions in integers.~~~~~~~~~~~~~~~~~~~~~~~~~~~~~~~~~~~~~~~~~~~~~~~~~~~~~~~~~~~~~~~~~~~~~~~~~~~~~~~~~~~~~~~~~~~~~~~~~~~~$\Box$ 

\noindent
{\bf Acknowledgement.} I thank M.A. Prasad for many
fruitful discussions and valuable suggestions.

\vspace{0.3cm}
\noindent
{\bf REFERENCES}
\begin{enumerate}

\item L.J. Mordell {\it Diophantine Equations},
ACADEMIC PRESS London and New York, 1969, Chapter 24.

\item Paulo Ribenboim {\it Fermat's Last Theorem for Amateurs},
Springer, 1999, Chapter VIII.

\end{enumerate}

\end{document}